\documentclass[11pt,a4paper]{article}
\usepackage[utf8]{inputenc}
\usepackage[english]{babel}
\usepackage{amsmath}
\usepackage{amsfonts}
\usepackage{amssymb}
\usepackage{makeidx}
\usepackage{graphicx}
\usepackage{lmodern}
\usepackage{kpfonts}
\usepackage{bbm}
\usepackage{import}
\usepackage{graphicx}
\usepackage{appendix}
\usepackage{comment}
\usepackage{appendix}
\usepackage{array}
\usepackage{enumitem}
\usepackage{authblk}
\usepackage[dvipsnames]{xcolor}
\usepackage[left=1.5cm,right=1.5cm,top=1.5cm,bottom=1.5cm]{geometry}
\newtheorem{theorem}{Theorem}
\newtheorem{definition}{Definition}
\newtheorem{proof}{Proof}
\newtheorem{proposition}{Proposition}

\newtheorem{remark}{Remark}

\newcommand{\beq}{\begin{eqnarray}}
\newcommand{\eeq}{\end{eqnarray}}

\newcommand{\bpr}{\begin{proof}}
\newcommand{\epr}{\end{proof}}

\newcommand{\ba}{\begin{align*}}
\newcommand{\ea}{\end{align*}}

\newcommand{\be}{\begin{equation}}
\newcommand{\ee}{\end{equation}}

\newcommand{\R}{\mathbb R}
\newcommand{\nn}{\nonumber}

\allowdisplaybreaks
\title{Duality for a boundary driven asymmetric model of energy transport}
\author[1]{Gioia Carinci}
\author[1,2]{Francesco Casini}
\author[1]{Chiara Franceschini}
\affil[1]{\small{Università di Modena e Reggio Emilia, FIM, Modena, Italy}}
\affil[1]{\small{Università di Parma , SMFI, Parma, Italy}}

\begin{document}
	
	\maketitle



{\begin{abstract}{We study the Asymmetric Brownian Energy, a model of heat conduction  defined on the one-dimensional finite lattice with  open boundaries. The system is shown to be dual to the Symmetric inclusion process with absorbing boundaries. The proof relies on  a non-local map  transformation procedure relating the model  to its symmetric version. As an application, we show how the duality relation can be used  to analytically compute suitable exponential moments  with respect to the stationary measure.}
\end{abstract}}

\section{Introduction}

The Asymmetric Brownian Energy Processes (ABEP) is an interacting diffusion system    describing an asymmetric energy exchange between the sites of a lattice.
Its symmetric version (BEP) was originally introduced in \cite{GKRV} where its symmetries and duality properties where unveiled. These are related to the intrinsic algebraic structure of the infinitesimal generator that can be written in terms of a continuous representation of the non-compact $\mathfrak{su}(1,1)$ Lie algebra.

In \cite{GKRV} the BEP in the closed system (i.e. in absence of external reservoirs) was proven to be dual to the  symmetric inclusion process (SIP). This is an interacting particle system modelling particles moving on a lattice with an attractive interaction. The reason behind the above mentioned duality relation lies is in the   $\mathfrak{su}(1,1)$ algebraic structure shared by the two processes.  BEP and  SIP are indeed two elements of a broader class of models all related to  the $\mathfrak{su}(1,1)$ Lie algebra and including also other  notable models. One of these is the Kipnis-Marchoro-Presutti model (KMP) \cite{KMP} where the total energy is instantaneously redistributed among sites and that can  be recovered as an instantaneous thermalization limit of the BEP. Another model inherently related to the BEP is the Wright-Fisher diffusion \cite{Ethe} that is the prototype model of mathematical population genetics. Duality between the Wright-fisher  and the Moran model can be seen as a particular instance of duality between BEP and SIP (see e.g. \cite{CGGR1}).

 In \cite{CGGR} the analysis was extended to the non-equilibrium situation in which the system  is put in contact with two external heat reservoirs imposing two different temperatures $T_\ell\neq T_r$ at the endpoints of the bulk. The corresponding process is also called   BEP with open boundaries and has been shown to be dual to the SIP with absorbing boundaries. 

Not long ago two new models belonging to this class have been introduced in \cite{Frassek1} via integrable non-compact spin chains and their duality relation shown in \cite{Frassek2, FFG}. For these new models formulae for the non-equilibrium steady states were obtained in \cite{Frassek3, FFG}. A characterization of these measures as mixture products of inhomogeneous distributions has been revealed in \cite{mixture1, mixture2}, however for an asymmetric dynamics this characterization is still an open problem.

The asymmetric version of the model we study (ABEP) was first introduced in \cite{CGSR2} in a closed boundary setting. This emerged as a scaling limit of the  ASIP (an asymmetric version of the inclusion process) in a particular regime of weak asymmetry. In the same work, an alternative construction was proposed for the ABEP, that was shown to be obtainable from BEP, via a non-local transformation $g$ depending on the asymmetry parameter. A duality relation between ABEP and SIP was then deduced in \cite{CGSR2}  as a consequence, independently, of the two above mentioned constructions. The duality function does not have a standard product structure (as is usually the case in the symmetric context) but a nested product structure related to the non-local map $g$.
This property is a first instance of a duality relation between   a genuinely non-equilibrium asymmetric system (in
the sense that it has a non-zero average current) and  a symmetric process. This link is made possible by the fact that the dependence on the asymmetry parameter is retained in the duality function, through the map $g$.

Here we extend the analysis  to the system with open boundaries. In this context  the problem becomes the definition of reservoirs itself. The aim is indeed to impose external temperatures $T_\ell\neq T_r$ in such a way as not to alter the condition of existence of a duality relation with the SIP with absorbing boundaries.
Our strategy does not directly rely  on algebraic considerations on the Markov generator but rather on the  link between ABEP and BEP via the non-local map $g$.
This transformation procedure allows us to construct reservoirs of the correct form. These turn out to act in a  non-standard way. The left reservoir acts only on the left endpoint of the lattice, but its action takes into account the total energy of the system. The right reservoir, instead affects all the sites of the lattice.
As a result of this construction we prove a duality relation with  the SIP with  absorbing boundary by means of two different duality functions. The first one is in a so-said {\em classical} form whereas the second one is in terms of {\em generalized Laguerre polynomials}.

As far as we know, duality in the presence of an asymmetry together with open boundary condition is still a quite challenging outcome as the classical techniques relying on algebraic considerations do not work. This is due to the fact that   the quantum group symmetry needed to construct the duality relation is broken. Results are mainly available for the case of asymmetric simple exclusion process (ASEP). The first attempt is due to Okhubo in \cite{Ohkubo} where a dual operator has been obtained; however it could not be directly interpreted as a transition matrix for a stochastic process. 
We mention \cite{Jeff} where the author generalizes the self-duality of the asymmetric simple exclusion process with an open boundary condition at the left boundary and a closed right boundary. More recent results include \cite{Guillaume} where a duality relation between an half-line open ASEP and a sub-Markov process where particles perform an asymmetric exclusion dynamics in the bulk and are killed at the boundary is proven. In \cite{Gunter, Gunter2} it is shown a reverse duality relation for an open ASEP with open boundary and a shock ASEP with reflecting boundary.

%
%

\vskip.2cm
\noindent
{The rest of the paper is organized as follows: in Section \ref{sec2} we introduce the model of interest, i.e.  the ABEP with open boundaries, and show how it can be obtained from its symmetric version via a non-local map transformation. At the end of the section 
we state some general results which allow to infer properties for a process that can be obtained from another process via a map transformation. 
In the subsequent two sections  these general properties are then specialized to gather information for ABEP starting from known results for BEP: Section \ref{revmes} is devoted to the study of the case $T_\ell=T_r$ in which the system is proven to be reversible, and the reversible measure is computed; in  Section \ref{dualitysection} instead we discuss duality relations. We end with Section \ref{sec6} where we use the duality results to gather some information on the stationary measure in the general case. In particular we compute what we call the  one-point and two-point  stationary exponential correlations of the partial energies.}

\section{The model} \label{sec2}
The Brownian Energy Process BEP$(\alpha)$ is an interacting
diffusion system of continuous spins placed on the sites of a   lattice $V$, $\alpha $ is a positive parameter tuning the intensity of the interaction.  We consider its asymmetric version ABEP$(\sigma,\alpha)$, $\sigma>0$ the asymmetry parameter,  that can be defined when the lattice is one-dimensional $V  = \lbrace 1, \ldots, N \rbrace $ and the interaction is  nearest-neighbor.  To each site of the lattice $i\in V$ is associated an energy $x_i\ge 0$. We denote by $x = (x_1, \ldots, x_N) \in\mathbb R_+^N$ the vector collecting all energies and we call  $\Omega  :=\mathbb R_+^N$ the state space of the system.
When the system is {\em closed}, or, in other words, in absence of external reservoirs, the dynamics conserves the total energy of the system $E(x):=\sum_{i\in V} x_i$.

In this paper we consider the {\em open system}, i.e. we put the {\em bulk} lattice $V$ in contact with two external reservoirs placed at artificial sites   $V^{{\rm res}} =\left\lbrace 0, N+1 \right\rbrace $. Each reservoir  $j \in V^{{\rm res}} $ can be interpreted as a thermal bath characterized by its own fixed temperature $T_j\ge 0$,  that is attached to the bulk $V$ only through the boundary sites 1 and $N$. The action of the reservoirs induces an energy exchange between the bulk lattice and the exterior,  that destroys the total energy conservation. For simplicity  we will also denote by $T_\ell:=T_0$ the temperature of the left reservoir and by $T_r:=T_{N+1}$ the temperature of the right reservoir.

\vskip.2cm
\noindent
In order to define the model, we need to define two crucial quantities the  partial energies $E_i (x)$, $i\in V$, and  the non-local map $g$.

\begin{definition} \label{mapg}
We define the map $g: \Omega  \rightarrow \Omega  $ via 
\begin{equation}
g(x) = \left(  g_i \left( x \right) \right)_{i \in V }  \qquad \text{with} \qquad g_i(x) := \dfrac{e^{-\sigma    E_{i+1}(x)} - e^{-\sigma    E_i(x)}}{\sigma   } 
\end{equation} 
where $E_i (x)$  denotes the energy of the system at the right of site $i \in V $, i.e.
\begin{equation}
E_i (x) = \sum_{l  = i}^{N} x_{l } \quad \text{for} \quad i=1,\ldots,N \qquad \text{with the convention} \quad E_{N+1}(x) = 0 \;.
\end{equation} 
\end{definition}
Notice that the total energy $E(x)$ coincides with the first component $E_1(x)$ of the vector of partial energies.

\vskip.2cm
\noindent
The  stochastic evolution of the collection of energies of the system is governed by a Markov process $\lbrace x(t), t\geq 0 \rbrace$ that we will define by giving its  infinitesimal generator $\mathcal{L}^{{\rm  ABEP}} $. This acts on smooth functions $f: \Omega\rightarrow \mathbb{R}$ and is given by the sum of three terms, one of them governing the interaction between bulk sites and the other two  modelling the action of  left and right reservoirs. We define
\begin{equation}\label{abep}
\mathcal{L}^{{\rm  ABEP}} = \mathcal{L}^{{\rm  ABEP}}_{{\rm left}} +\sum_{i=1}^{N-1} \mathcal{L}^{{\rm  ABEP}}_{i,i+1} + \mathcal{L}^{{\rm  ABEP}}_{{\rm right}}.
\end{equation}
where, for $i\in \{1, \ldots, N-1\}$, the action on smooth functions $f:\Omega \to \mathbb R$ is
\begin{align} \label{abep-bulk}
 [\mathcal{L}^{{\rm  ABEP}} _{i,i+1} f](x) &=    \dfrac{1}{2\sigma   ^{2}}  \left( 1 - e^{-\sigma    x_{i}}\right) \left( e^{\sigma    x_{i+1}} -1\right)  \left( \dfrac{\partial}{\partial x_{i+1}} - \dfrac{\partial}{\partial x_{i}} \right) ^{2} f(x) \\ &   + \dfrac{1}{\sigma   }  \biggl(  \left( 1 - e^{-\sigma    x_{i}}\right) \left( e^{\sigma    x_{i+1}} -1\right) + \alpha \left(2 - e^{-\sigma    x_i} - e^{\sigma    x_{i+1}} \right)  \biggr) \left( \dfrac{\partial}{\partial x_{i+1}} - \dfrac{\partial}{\partial x_{i}} \right) f(x)\nonumber
\end{align}
whereas
\begin{equation} \label{abep-res}
[\mathcal{L}^{{\rm  ABEP}}_{{\rm left}}f](x)= T_\ell \left(e^{\sigma    E(x) }  \left( \alpha -1 + e^{\sigma    x_1} \right)  \dfrac{\partial}{\partial x_{1} } + \dfrac{e^{\sigma    E(x)}}{\sigma   }\left(e^{\sigma    x_{1}} - 1 \right) \dfrac{\partial^{2}}{\partial x_{1}^{2}} \right) f(x)-  \dfrac{e^{\sigma    x_{1}} -1}{\sigma   } \dfrac{\partial}{\partial x_{1}}f(x)
\end{equation}
and
\begin{align} \label{abep-res-dx}
[\mathcal{L}^{{\rm  ABEP}}_{{\rm right}}f](x)&= \left(\alpha T_r- \dfrac{1-e^{-\sigma    x_N}}{\sigma   } \right) \sum_{l =1}^{N} e^{\sigma    E_{l }(x)} \left( \partial_{x_l } - \partial_{x_{l -1}} \right)  f(x) 
  + \\ & T_r \dfrac{1-e^{-\sigma    x_N}}{\sigma   }  \sum_{l , j =1}^{N} e^{\sigma    \left( E_{l } (x) + E_{j} (x) \right)} \left( \partial_{x_l } - \partial_{x_{l -1}}  \right)  \left( \partial_{x_j} - \partial_{x_{j-1}}  \right)+ T_r \left( 1-e^{-\sigma    x_N}  \right) \sum_{l  = 1}^{N} e^{2\sigma    E_{l }(x) }  \left( \partial_{x_l } - \partial_{x_{l -1}}  \right)  f(x) \ .
  \nonumber
\end{align}


The action of reservoirs is non-local in two different ways. The left reservoir acts only on the left boundary site 1, but its action takes in account the total energy $E(x)$ that is not an invariant of the dynamics. The right reservoir, instead,  affects the whole chain.

\section{From BEP to ABEP}

The BEP$(\alpha)$ on $V$ is the symmetric version  of the ABEP$(\sigma, \alpha)$  obtained  in the limit as $\sigma \to 0$. As in the previous section we consider the system with nearest-neighbor interaction in contact with two boundary reservoirs kept at temperature $T_\ell$ and $T_r$. We denote by $\{z(t), \ t \ge 0\}$ the Brownian Energy  process on the space state $\Omega= \mathbb R_+^N$ describing the evolution of the vectors $z:=(z_1, \ldots, z_N)$ of single-site energies. The infinitesimal generator, acting on smooth functions $f:\Omega \to \mathbb R$, is defined as follows
    \begin{equation}\label{bep}
      \mathcal{L}^{{\rm BEP}} =  \mathcal{L}^{{\rm  BEP}}_{{\rm left}}+ \sum_{i=1}^{N-1} \mathcal{L}_{i,i+1}^{{\rm BEP}} +
       \mathcal{L}^{{\rm  BEP}}_{{\rm right}}
    \end{equation}
    where, for $i\in \{1, \ldots, N-1\}$, 
    \begin{equation} \label{bulkbep}
		\mathcal{L}_{i,i+1}^{{\rm BEP}} f (z) = \left[  z_{i}z_{i+1}\left(\partial_{z_{i+1}}-\partial_{z_{i}}\right)^{2}-\alpha(z_{i}-z_{i+1})\left(\partial_{z_{i+1}}-\partial_{z_{i}}\right) \right] f(z)
	\end{equation}
	whereas 
    \begin{equation}
        \mathcal{L}^{{\rm  BEP}}_{{\rm left}} f (z) = \left[  T_\ell \left( \alpha  \dfrac{\partial}{\partial_{z_{1}}} + z_{1} \dfrac{\partial^{2}}{\partial_{z_{1}}^{2}} \right)  -  z_{1} \dfrac{\partial}{\partial_{z_{1}}} \right] f (z) 
        \end{equation}
        and
        \begin{equation}\label{end}
          \mathcal{L}^{{\rm  BEP}}_{{\rm right}} f (z)  = \left[  T_r \left( \alpha  \dfrac{\partial}{\partial_{z_{N}}} + z_{N} \dfrac{\partial^{2}}{\partial_{z_{N}}^{2}} \right)  -  z_{N} \dfrac{\partial}{\partial_{z_{N}}} \right] f (z) \ .
    \end{equation}
    The  latter terms give the action of    left and  right reservoirs that are attached, respectively, to site $1$ and site $N$.

   \vskip.3cm
   \noindent 
    It can be easily checked that $ \mathcal{L}^{{\rm BEP}}$ is recovered from $ \mathcal{L}^{{\rm ABEP}}$ by suitably taking the limit as $\sigma \to 0$. On the other hand $ \mathcal{L}^{{\rm ABEP}}$ can be constructed from $ \mathcal{L}^{{\rm BEP}}$ by acting with the non-local map $g$ introduced in Definition \ref{mapg}. This claim has been proven in \cite{CGSR2} for the closed system. Below we  show that such a construction can be extended to the reservoir terms of the generator.

\begin{theorem}[From BEP to ABEP]
Let $g$ be the map  in Definition \ref{mapg}, then for all   $f \in {\cal D}({\cal L}^{{\rm BEP}})$ we have
\begin{equation}\label{connect}
	 \mathcal{L}^{{\rm ABEP}} (f \circ g) =\left[ \mathcal{L}^{{\rm BEP}} f \right] \circ g \  .
		\end{equation}
\end{theorem}
\textbf{Proof:} Throughout this proof we will use the alternative notation for the reservoir terms of the generators
$\mathcal{L}^{{\rm  BEP}}_{0,1}:=\mathcal{L}^{{\rm  BEP}}_{{\rm left}}$ and $ \mathcal{L}^{{\rm  BEP}} _{N,N+1}:=\mathcal{L}^{{\rm  BEP}}_{{\rm right}}$, respectively, $\mathcal{L}^{{\rm  ABEP}}_{0,1}:=\mathcal{L}^{{\rm  ABEP}}_{{\rm left}}$ and $ \mathcal{L}^{{\rm  ABEP}} _{N,N+1}:=\mathcal{L}^{{\rm  ABEP}}_{{\rm right}}$ and write, for $(i,j)\in V \times V^{{\rm res}}$, 
\begin{equation}\label{three}
		\mathcal{L}_{i,j}^{{\rm BEP}} f (z) =  \left[  T_{j} \left( \alpha  \dfrac{\partial}{\partial_{z_{i}}} + z_{i} \dfrac{\partial^{2}}{\partial_{z_{i}}^{2}} \right)  -  z_{i} \dfrac{\partial}{\partial_{z_{i}}} \right] f(z)
	\end{equation}
	and
	\begin{align} \label{abep-res-general}
[ \mathcal{L}^{{\rm  ABEP}} _{i,j} f](x)&= \left(\alpha T_j - g_i(x) \right) \left[ \sum_{l =1}^{i-1} e^{\sigma   E_{l }(x)} \left( 1-e^{-\sigma    x_l } \right) \partial_{x_l } + e^{\sigma    E_i (x)}  \partial_{x_i}  \right] f(x)+ \\
\nonumber  &
  T_j g_i(x) \left[   \sum_{l , j =1}^{i-1} e^{\sigma    E_{l } (x)} \left( 1- e^{-\sigma    x_l } \right)  \, e^{\sigma    E_{j} (x)} \left( 1- e^{-\sigma    x_j} \right) \partial^{2}_{x_l  x_j}  + e^{2\sigma   E_i (x)} \partial^{2}_{x_i}  \right. \\ 
 \nonumber & \left.  
 2 \sum_{l =1}^{i-1} e^{\sigma    E_{l }(x)} \left(  1 - e^{-\sigma    x_l }\right)e^{\sigma    E_i (x)} \partial^{2}_{x_l  x_i} + \sum_{l  =1}^{i-1} \sigma    e^{2\sigma    E_{l } (x) } \left( 1- e^{-2\sigma    x_l } \right) \partial_{x_l } + \sigma    e^{2\sigma    E_i (x)} \partial_{x_i} 
 \right] f(x) \;. 
\end{align}
In this way we can resume the proof of the theorem in the following two steps:
\begin{enumerate}
    \item for all $i \in V $, $x \in \Omega$,
	\begin{equation}
	\left[ \mathcal{L}_{i,i+1}^{{\rm BEP}} f \right] (g(x)) = 	\left[ \mathcal{L}_{i,i+1}^{{\rm ABEP}} f \circ g  \right] (x)
		\end{equation}
		
  \item for all $(i,j)\in V \times V^{{\rm res}}$, $x \in \Omega$,
 	\begin{equation}\label{compo}
	\left[ \mathcal{L}_{i,j}^{{\rm BEP}} f \right] (g(x)) = 	\left[ \mathcal{L}_{i,j}^{{\rm ABEP}} f \circ g  \right] (x)
		\end{equation}
\end{enumerate}
Step 1 has been proven in Theorem 3.4 of  \cite{CGSR2}.
It remains to prove Step 2.
Recalling the definition of $g$:
\begin{align*}
	g:\Omega  \;&\rightarrow\; \Omega \\
	x\;&\rightarrow\;g(x)=(g_{i}(x))_{i\in V}, \qquad \text{with} \qquad g_i(x)=\frac{e^{-\sigma    E_{i+1}(x)}-e^{-\sigma    E_{i}(x)}}{\sigma   }
\end{align*}
with the convention $E_{N+1}(x)=0$ and $E_1(x)=E(x)$.

\vskip.2cm
\noindent
Notice that the map $g$ is not full range, i.e. $g(\Omega )\neq \Omega $, indeed

\begin{equation}
E(g(x))= \frac 1 {\sigma} \left(1-e^{-2\sigma E(x)}\right)\le \frac1{\sigma}
\end{equation}
so that in particular $g(\Omega)\subseteq \{x\in \Omega : E(x)\le 1/\sigma\}$.
Moreover  $g$ is a bijection from $\Omega$
to $g(\Omega)$. Indeed,  denoting by $g^{\rm{inv}}: g(\Omega)  \rightarrow \Omega$ the  inverse transform of $g$.  In other words, if $z=g(x)\in g(\Omega)$, then   $x=g^{\rm{inv}}(z)$
with $i$th component being
\be
g^{\rm{inv}}_i(z)=\frac 1 {\sigma   } \, \ln \left\{ \frac{1-\sigma    E_{i+1}(z)}{1-\sigma    E_{i}(z)}\right\}
\ee

Let $F:= f \circ g$, or, equivalently, $f=F \circ g^{\rm{inv}}$ namely
$
F(x)=f(g(x)) $ for $x\in \Omega$ and $ f(z)=F(g^{\rm{inv}}(z)) $ for $z\in g(\Omega)$, therefore, in order to prove \eqref{compo}, it is sufficient to show that, for all $x \in \Omega$,
	\begin{equation}\label{compo1}
	\left[ \mathcal{L}_{i,j}^{{\rm BEP}} (F \circ g^{\rm{inv}} )\right] (g(x)) = 	\left[ \mathcal{L}_{i,j}^{{\rm ABEP}} F \right] (x)
		\end{equation}
At this aim we compute the first and second  derivatives of $f=F \circ g^{\rm{inv}}$. We have 
\beq \label{Df}
\frac{\partial f}{\partial z_k}(z)= \sum_{l \in V} \frac{\partial F}{\partial x_l }(g^{\rm{inv}}(z)) \cdot \frac{\partial g^{\rm{inv}}_l  }{\partial z_k}(z) \qquad \text{for all} \qquad k \in V
\eeq
and
\beq \label{DDf}
\frac{\partial^2 f}{\partial^2 z_k z_m}(z)= 
\sum_{l ,j\in V} \frac{\partial^2 F}{\partial^2 x_l  x_j}(g^{\rm{inv}}(z)) \cdot \frac{\partial g_{l}^{\rm{inv}} }{\partial z_k}(z) \cdot \frac{\partial g_{j}^{\rm{inv}}}{\partial z_m}(z) 
+ \sum_{l \in V} \frac{\partial F}{\partial x_l }(g^{\rm{inv}}(z)) \cdot \frac{\partial^2 g^{\rm{inv}}_{l} }{\partial^2 z_k z_m}(z) \qquad \text{for all} \qquad k,m \in V \;.
\eeq
We now compute all the first and second derivatives of all the components of the inverse function $g^{\rm{inv}}$, we obtain
\begin{equation}\label{firstDerivativeh}
\frac{\partial g^{\rm{inv}}_l }{\partial z_k}(z)=\begin{cases}
	0\quad &\text{if}\;k<l \\
	\frac{1}{1-\sigma    E_{l }(z)}\quad &\text{if}\;k=l \\
	\frac{\sigma    z_{l }}{(1-\sigma    E_{l }(z)) (1-\sigma    E_{l  +1}(z))}\quad &\text{if}\;k>l 
\end{cases}
\end{equation}
and,  for $m\leq k$ (it is symmetric in $k$ and $m$),
\begin{equation}\label{secondDerivativeh}
	\frac{\partial^{2} g^{\rm{inv}}_{l }}{\partial z_{m}\partial z_{k}}(z)= \begin{cases}
		0\quad &\text{if}\; m<l \\
		\frac{\sigma   }{\left(1-\sigma    E_{l }(z)\right)^{2}}\quad &l =m\leq k\\
		\frac{z_{l }\sigma ^{2}(2-\sigma  z_{l }-2\sigma    E_{l +1}(z))}{\left[(1-\sigma    E_{l }(z)) (1-\sigma    E_{l  +1}(z))\right]^{2}} \quad &\text{if}\;m>l 
		\end{cases}
\end{equation}
These derivatives simplify observing that, thanks to telescopicity of the sum,  
%
\begin{equation}\label{energyZenergyX}
E_{l }(z)=\sum_{i=l }^{N}z_{i}=\sum_{i=l }^{N}g_{i}(x)=\frac{1}{\sigma   }\left(1-e^{-\sigma    E_{l }(x)}\right) \;.
\end{equation}
And then, using \eqref{energyZenergyX} we can simplify the expressions for  the derivatives as  follows:
\begin{equation}\label{Dh}
\frac{\partial g^{\rm{inv}}_{l }}{\partial z_{k}}(z)=
\begin{cases}
0\quad &\text{if}\;k<l \\
	e^{\sigma    E_{l }(x)}\quad &\text{if}\;k=l \\
	e^{\sigma    E_{l }(x)} \left( 1 - e^{-\sigma    x_{l }} \right) \quad &\text{if}\;k>l 
\end{cases}
\end{equation}
and 
\begin{equation}\label{DDh}
\frac{\partial^{2} g^{\rm{inv}}_{l }}{\partial z_{m}\partial z_{k}}(z)= \begin{cases}
		0\quad &\text{if}\; m<l \\
		\sigma    e^{2\sigma    E_{l }(x)}\quad &l =m\leq k\\
		\sigma    e^{2\sigma    E_{l }(x)} \left(1-e^{-2\sigma    x_{l }}\right)  \quad &\text{if}\;m>l  \ .
		\end{cases}
\end{equation}
Then  by substituting the  expressions \eqref{Dh}
and \eqref{DDh} into equations \eqref{Df} and \eqref{DDf} we obtain  explicit expressions for the first and second  derivatives of $f=F \circ g^{\rm{inv}}$. 
Finally, the identity  \eqref{compo1} follows by replacing these expressions into the BEP$(\alpha)$ boundary generators given in  \eqref{three}.
\begin{flushright}
    $\square$
\end{flushright}

\subsection{Some general definitions and properties}\label{sec3}

The construction of ABEP$(\sigma,\alpha)$ as a non-local transformation of BEP$(\alpha)$, allows to derive several fundamental properties of the asymmetric process, such as  duality properties or the structure of the stationary measure. These are by starting from the analogous properties of the symmetric process and projecting them via the map $g$. 
Having this goal in mind, in this section we prove some general results  relating two Markov processes that are connected  via a map transformation. 
\vskip.2cm
\noindent
We start by recalling  the definition duality in terms of infinitesimal generators of two Markov processes. We will denote by ${\cal D}({\cal L})$ the domain of ${\cal L}$.
\vskip.2cm
\noindent

\begin{definition}[Generator duality]
\label{gen-daulity-abcd}
Let ${\cal L}$ and $L$ be the infinitesimal generators of two Markov  processes $\{X(t):t\geq 0\}$  and $\{Y(t): t\geq 0\}$ defined, respectively, on the state spaces $\Omega$ and $\Omega^{{\rm dual}}$. Let $D: \Omega\times {\Omega}^{{\rm dual}}\to\R$ be a measurable function, such that
$D(y, \cdot) \in {\cal D}({\cal L})$ and $D(\cdot, x)\in {\cal D}(L)$.
We then say that $D$ is a duality function for generator duality between the processes
$\{X(t):t\geq 0\}$ and $\{Y(t): t\geq 0\}$ if
for all $x\in \Omega, y\in {\Omega}^{{\rm dual}}$, we have
\be
\label{gendualfirstdef}
\left({\cal L} D(\cdot, x)\right) (y)=\left(L D(y, \cdot) \right)(x)
\ee
\end{definition}

 \vskip.2cm
\noindent
In the next theorem we will see that if a stationary measure,  a reversible measure or a duality function are known for one of the a processes, then  the corresponding object  can be computed for a process obtained from the original one via a transformation. 


\begin{theorem}\label{General_result}
Let $g$ be a map $g:\Omega \to \Omega$, with $\Omega \subseteq \mathbb R_+^N$ and let ${\cal L}$ and $\widehat{\cal L}$ be the infinitesimal generators of two Markov processes on the state spaces, respectively $\Omega$ and $\widehat \Omega:=g(\Omega)$. Suppose that $\forall f \in {\mathcal D}({\cal L})$ it holds that $f \circ g \in  {\mathcal D}(\widehat {\cal L})$ and	\begin{equation}\label{compo}
		\widehat{\cal L} (f\circ g)=\left({\cal L}f\right)\circ g
	\end{equation}
 then we have the following properties.
\label{main}
\begin{itemize}
\item[i)]  Let $\mu$ be a measure on $\Omega$ absolutely continuous w.r.t. Lebesgue. Let ${\cal J}$ be the Jacobian matrix of the map $g$. 
If $\mu$ is a  stationary (reversible) measure for ${\cal L}$  then 
\be\label{compomu}
\hat\mu :=  \left(\mu \cdot \rm{det}{\cal J}  \right)\circ g 
\ee
is a stationary (reversible) measure for 
$\widehat{\cal L}$. 
\item[ii)] Let $L$ be the infinitesimal generators a Markov processes on the state space $\Omega^{{\rm dual}}$. If ${\cal L}$ is dual to $L$ with duality function $D: \Omega \times \Omega^{{\rm dual}} \to \mathbb R$, then $\widehat{\cal L}$ is dual to $L$ with duality function $D: \widehat \Omega \times \Omega^{{\rm dual}} \to \mathbb R$
\begin{equation}
\widehat D(\cdot,\xi):= D(\cdot, \xi) \circ g \ ,\qquad \qquad \xi \in \Omega^{{\rm dual}} \ .
\end{equation}
\end{itemize}

\end{theorem}
\textbf{Proof:}
\vskip.1cm
\noindent
\begin{itemize}
\item[i)] Due to the  absolute continuity of $\mu$ we can write, with a slight abuse of notation, that $\mu(dx)= \mu(x) \ dx $. The stationarity condition for $\mu$ with respect to ${\cal L}$ then reads
\be
\int[{\cal L}f](z) \, \mu(z) \ dz  =0, \qquad \text{for all} \quad f \in {\cal D}({\cal L})
\ee
that, taking the change of variables $z=g(x)$, gives 
\be
\int [{\cal L}f](g(x)) \, \cdot \mu(g(x)) \cdot \rm{det} {\cal J}(g(x)) \, dx  =0, \qquad \text{for all} \quad f \in  {\cal D}({\cal L})
\ee
that, thanks to \eqref{compo} and \eqref{compomu}, is equivalent to
\be
\int [\widehat{\cal L}(f\circ g)](x) \, \cdot \hat\mu(x) \, dx  =0, \qquad \text{for all} \quad f \in  {\cal D}({\cal L}) \ .
\ee
Due to the fact that  $D(\widehat{\cal L})= \{F=f\circ g: \: f\in D({\cal L})\}$, the last identity can be rewritten as
\be
\int [\widehat{\cal L} F](x) \, \cdot \hat\mu(x) \, dx  =0, \qquad \text{for all} \quad F\in  {\cal D}(\widehat{\cal L})
\ee
that is the stationary condition of $\hat \mu$ with respect to $\widehat{\cal L}$. The statement regarding reversible measures can be proven in an analogous way.

    \item[ii)] To prove the second statement we use the duality relation between ${\cal L}$ and $L$ and take the composition of the duality function (as a function of the variable $x$) with the function $g$. For $x\in \Omega$ and $\xi \in \Omega^{{\rm dual}}$, we have
\begin{eqnarray}
\left[\widehat {\cal L}\widehat D(\cdot,\xi)\right](x)
&=&\left[\widehat{\cal L}(D(\cdot,\xi)\circ g)\right](x) \\
& = &\left[{\cal L}D(\cdot,\xi)\right](g(x)) \nn \\
&=&\left[LD(g(x),\cdot)\right](\xi)\nn \\
& =&\left[L \widehat D(x,\cdot)\right](\xi)\ .
\end{eqnarray}
 This concludes the proof of the second item.
\end{itemize}
\begin{flushright}
    $\square$
\end{flushright}
In the next two sections we specialize the argument of the above theorem for our model of interest. In Section \ref{revmes} we focus on the cas in which the external reservoirs impose the same temperatures (i.e. when $T_\ell=T_r=T$). We prove that in this situation ABEP$(\sigma,\alpha)$ is reversible and we find the reversible measure. In Section \ref{dualitysection} we find two duality relations for ABEP$(\sigma,\alpha)$.

\section{Equal temperature reservoirs} \label{revmes}
In this section use item i) of Theorem \ref{General_result} to withdraw some conclusions concerning the case in which   the two reservoirs have the same temperature. The idea is to import this property from the   reversibility of the corresponding symmetric process. From Section 3 of \cite{CGGR}  we know indeed that, in absence of reservoirs, the BEP$(\alpha)$ is reversible. In particular it admits a one-parameter family of reversible probability measures $\mu_T$, $T\ge 0$, that are products of   Gamma distributions of shape parameters $\alpha $ and scale parameter $T$, i.e. $\mu_{T}^{{\rm BEP}} (z) \ dz$ with
\be
\mu_{T}^{{\rm BEP}} (z)= \prod_{i=1}^{N}\frac{e^{-z_i /T}z_i^{(\alpha -1)}}{\Gamma(\alpha )T^{\alpha }}
\ee
When the process is in contact with two reservoirs kept at equal temperatures, $T_\ell=T_r=T$, the process remains reversible, admitting $\mu_T$ as the unique stationary probability measure. In the following Theorem we extend the statement to the asymmetric process, for which we prove the existence of a unique reversible probability  measure that is in the form of a product measure times a function of the total energy of the system $E(x)$.

\begin{theorem}[Reversible measure for ABEP]\label{rev}
The ABEP$(\sigma,\alpha)$ with equal reservoir temperatures $T_\ell = T_r = T $ is reversible with respect to the unique stationary probability measure $\mu_{T}^{{\rm ABEP}} (x) \ dx$, with
   \begin{equation}\label{muabep}
	\mu_{T}^{{\rm ABEP}} (x)= \exp\left\{ \frac{e^{-\sigma    E(x)} -1}{\sigma    T}  \right\}\cdot \prod_{i=1}^{N} \dfrac{(1-e^{-\sigma    x_i})^{(\alpha -1)} }{ \Gamma(\alpha ) \sigma^{\alpha -1} T^{\alpha }} \ e^{-\sigma  x_i(\alpha(i  -1)+1)} 
\end{equation} 
\end{theorem} 

\textbf{Proof:}
We want to use item i) of Theorem \ref{General_result}. To this aim it is enough to compute 
$(\mu_{T}^{{\rm BEP}}\circ g)(x)$. Indeed, 
\begin{align*}
	\mu_{T}^{{\rm ABEP}} (x) &= (\mu_{T}^{{\rm BEP}}\circ g)(x) = \prod_{i=1}^{N}\mu^{{\rm BEP}}_{T}(g_i(x))=\prod_{i=1}^{N}\frac{e^{-g_{i}(x) /T}(g_{i}(x))^{(\alpha -1)}}{\Gamma(\alpha )T^{\alpha }} \: {\rm det} \mathcal{J}(g(x))  \\& =
\prod_{i=1}^{N} \dfrac{1}{\Gamma(\alpha ) T^{\alpha }} \cdot \exp\left\{-\frac{e^{-\sigma    E_{i+1}(x)}-e^{-\sigma    E_{i}(x)}}{\sigma    T}\right\} (1-e^{-\sigma    x_i})^{(\alpha -1)} \ \dfrac{e^{-\sigma    (\alpha -1) E_{i+1}(x)}}{\sigma   ^{(\alpha -1)}} \:e^{-\sigma    E_i(x)}   \\& =
 \sigma   \cdot  { e^{\sigma  \alpha E(x)} e^{(T-\sigma   )E(x)}
\exp\left\{ \frac{e^{-\sigma    E(x)} -1}{\sigma    T}  \right\} }\cdot
\prod_{i=1}^{N} \frac{(1-e^{-\sigma    x_i})^{(\alpha -1)} e^{-(\sigma  \alpha i + T)x_i}  }{ (\sigma    T)^{\alpha }\Gamma(\alpha )} \\& =
\exp\left\{ \frac{e^{-\sigma    E(x)} -1}{\sigma    T}  \right\} \cdot \prod_{i=1}^{N} \dfrac{(1-e^{-\sigma    x_i})^{(\alpha -1)} }{\sigma ^{\alpha -1} T^{\alpha } \Gamma(\alpha )} e^{-(\sigma  \alpha i + T)x_i} e^{(\sigma  \alpha + T -\sigma    )x_i} 
\end{align*}   
here, by calling again $z=g(x)$, $J$ is the Jacobian $N\times N$ matrix given by 
\begin{equation}
	J(z)=\left(\frac{\partial g_{l}^{\rm{inv}}}{\partial z_{k}}\right)_{k\in\{1,\ldots,N\},\,l\in\{1,\ldots,N\}}
\end{equation}
where the partial derivative are computed in \eqref{Dh}. Therefore, \eqref{muabep} follows. 
\begin{flushright}
    $\square$
\end{flushright}

\begin{remark}
In Theorem 3.3 of \cite{CGSR2}) a family of reversible measures has been found for ABEP$(\sigma,\alpha)$ with closed
 boundary. This family is labeled by the temperature $T$. The measure corresponding to the temperature $T$ (eq. (3.15 - 3.16) of \cite{CGSR2}) does not match with $\mu_{T}^{{\rm ABEP}}$ found in \eqref{muabep}. Indeed it differs from it only for the factor in front of the product in \eqref{muabep} that is a function of the total energy $E(x)$. This is due to the fact that, in absence of reservoirs, the total energy is an invariant of the dynamics, and then this term becomes a constant that simplifies with the normalizing factor of the probability measure.  In the presence of two thermal reservoirs instead, even in the case of equal temperatures $T_\ell=T_r=T$, the system does not conserve the total energy anymore, and the initial factor in \eqref{muabep} can not be neglected anymore.
\end{remark}

\section{Duality results}\label{dualitysection}

When $T_\ell \neq T_r$  reversibility is lost. Nevertheless there exists a unique  stationary measure   depending on both temperatures $T_\ell$ and $T_r$. However a full characterization of such a  measure is a difficult and still  open problem, even for the symmetric case. A tool that has proven to be of great help in the study of the properties of the  stationary measure is duality. We will return to the study of   steady state in Section \ref{sec6}. In the next section we prove two duality relations between the Asymmetric Brownian Energy process and the Symmetric Inclusion process with absorbing boundaries.

\subsection{Duality between ABEP and SIP}
The Symmetric Inclusion Process is a system of interacting particles moving in a lattice with attractive, nearest-neighbor interaction. It was originally introduced in \cite{GKR} as the attractive counterpart of the Simple Symmetric Exclusion Process. Each site can host for an unbounded number of particles, and then the state space of the inclusion process on the lattice $V=\{1, \ldots, N\}$ is $\mathbb N_0^N$. The attraction intensity is tuned by a parameter $\alpha>0$. Each particle may jump to its left or its right with rates proportional to the number of particles in the departure site and to the number of particles in the arrival site plus $\alpha$. We use the acronym SIP$(\alpha)$ for the Symmetric Inclusion process of parameter $\alpha$. Duality between BEP$(\alpha)$ and SIP$(\alpha)$ is well known in the literature. When the BEP system is put in contact with two external reservoirs a  duality relation still holds true. The dual process is still a system of inclusion particles inclusion, with the difference that the boundary conditions at the endpoints of the chain are no longer closed but absorbing.
 We give below the definition of the SIP$(\alpha)$ with absorbing boundaries. Notice that for this process the boundary sites $0$ and $N+1$ are no longer {\em artificial}, since their state is relevant in the dynamics. Configurations are then $N+2$-dimensional vectors that will be denoted by $\xi:=(\xi_0, \xi_1, \ldots, \xi_N, \xi_{N+1})$, $\xi_i$ being the number of particles at site $i$. The state space of this process is then the set $\Omega^{{\rm dual}}=\mathbb N_0^{V\cup V^{{\rm res}}}$, keeping in mind that, even if we keep the same notation  $V^{{\rm res}}$ for the set $\{0,N+1\}$, these sites in the dual process do no longer have  the meaning of reservoirs  but represent the absorbing sites. These leave eventually the bulk empty by absorbing all the particles.
\begin{definition} [SIP with absorbing boundaries]\label{def-sip}
We denote by $\{\xi(t), \ t \ge 0\}$ the SIP$(\alpha)$ on $V$ with absorbing boundaries $0$ and $N+1$, the Markov process on  $\Omega^{{\rm dual}}=\mathbb N_0^{V\cup V^{{\rm res}}}$ whose  infinitesimal generator acts on functions $f: \Omega^{{\rm dual}} \rightarrow \mathbb{R}$  and is defined as follows: 
\begin{equation}\label{sip}
{L}^{{\rm SIP-abs}}= {L}^{\rm{abs}}_{{\rm left}} + \sum_{i =1}^{N-1} {L}_{i,i+1}^{{\rm SIP}} + {L}^{\rm{abs}}_{{\rm right}}   \;,
\end{equation}
where, for all   and $i\in \{1, \ldots N-1\}$,
\begin{align} \label{sipbulk}
[{L}_{i,i+1}^{{\rm SIP}} \ f ]  (\xi )=  \sum_{i=1}^{N-1}  \xi_i (\alpha +\xi_{i+1}) \left[ f(\xi^{i,i+1})  - f(\xi)\right] + \xi_{i+1} (\alpha +\xi_{i}) \left[ f(\xi^{i+1,i})  - f(\xi)\right] 
\end{align}
and 
\begin{equation}
[{L}^{\rm{abs}}_{{\rm left}}f](\xi) := \xi_1 \left[ f(\xi^{1,0})  - f(\xi)\right] \qquad \text{and}\qquad [ {L}^{\rm{abs}}_{{\rm right}}  \ f ] (\xi )= \xi_N \left[ f(\xi^{N,N+1})  - f(\xi)\right] \; .
\end{equation}
\end{definition}

\vskip.2cm
\noindent
Besides being dual to the BEP, the Inclusion process with closed boundaries has been proved to be dual to the ABEP. This property has been proved in   \cite{CGSR2} and is  the first example of duality between an asymmetric system (i.e. bulk-driven) and a symmetric system (with zero current). This is made possible by the fact that the dependence on the asymmetry parameter $\sigma$ is  transferred to the duality function.  Here we generalize the result to the ABEP with reservoirs, that will be proven to be dual, again, to the Inclusion process with absorbing boundaries, exactly as its symmetric counterpart. This property will be proven using item 2 of Theorem \ref{General_result} and using the relation \eqref{General_result} that connects ABEP and BEP through the map $g$. We will prove two different duality relations between the same two processes. The first  duality relation is via the so-called {\em classical duality function} \cite{CGGR}, the second is in terms  of a duality function that is a product of  Laguerre polynomials, i.e. of the type {\em orthogonal polynomial duality function} \cite{simone}.

\subsubsection{Duality properties for the symmetric process.}

We start by  showing two duality relations between BEP$(\alpha)$ with reservoirs and SIP$(\alpha)$ with absorbing boundaries. The relations are given with respect to two different duality functions. The first duality relation is well-known, it is given in terms of the so-called {\em classical duality} and has been proven in \cite{CGGR}. The second result instead is given in terms of a duality function belonging to the class of {\em orthogonal polynomials dualities}, and more precisely it is related to the so-called {\em generalized Laguerre polynomials}.
Differently from the classical one, the orthogonal duality result for the open system is new, being available only for the closed system (see \cite{FG} for the proof).

\begin{theorem}[Duality between open BEP and  SIP with absorbing boundaries]\label{symduality}
The BEP$(\alpha)$ with an open  boundaries, with generator ${\cal L}^{{\rm BEP}}$ defined in \eqref{bep}-\eqref{end}, is dual to the SIP$(\alpha)$ with absorbing boundaries  defined in Definition \ref{def-sip} with respect to the following duality functions:
\begin{enumerate}
\item {\bf classical duality:} 
\begin{equation} \label{df}
D(z, \xi) = T_\ell ^{\xi_0} \cdot \prod_{i=1}^{N} \dfrac{\Gamma(\alpha )}{ \Gamma(\alpha  + \xi_i)} z_i^{\xi_i} \cdot  T_r^{\xi_{N+1}}, 
\end{equation}
\item {\bf orthogonal duality:}
\begin{equation}\label{odbepsip}
{\mathfrak D}_T(z,\xi) = \left(  T_\ell - T \right)^{\xi_0} \cdot \prod_{i=1}^{N} (-T)^{\xi_i}\cdot \mathstrut_1 F_1 \left( {\left. \genfrac{}{}{0pt}{} {-\xi_{i}} { \alpha }  \right\vert {\frac{z_{i}}{T}}} \right)\cdot
 (T_r-T)^{\xi_{N+1}}\;,
\end{equation}
for all $T >0$.
\end{enumerate}
\end{theorem}
\textbf{Proof:} 
For the proof of item 1 we refer to Theorem 4.1 of \cite{CGGR}.
In order to prove the second item, we have to show that
\begin{equation}
 	[\mathcal{L}^{{\rm BEP}}  {\mathfrak D}_T(\cdot, \xi) ](z) = [{L}^{{\rm SIP}} {\mathfrak D}_T(z, \cdot) ](\xi)
\end{equation}
Since both $\mathcal{L}^{{\rm BEP}} $ and ${L}^{{\rm SIP}}$ of a bulk term and two reservoir terms, it is sufficient to show that the duality relation for generators holds true term by term. The relation for the bulk terms of the generators has been proved in Section 4.2  of \cite{FG}, where it has been shown that, defining  $d(\zeta, k) =  (-T)^{k}\mathstrut_1 F_1 \left( {\left. \genfrac{}{}{0pt}{} {-k} { \alpha }  \right\vert {\frac{\zeta}{T}}} \right) $, for all $i\in \{1, \ldots, N-1\}$,
\begin{equation}
 [\mathcal{L}_{i,i+1}^{{\rm BEP}} d_T(\cdot, \xi_i) \cdot  d(\cdot, \xi_{i+1})] (z_i, z_{i+1}) = [{L}_{i,i+1}^{{\rm SIP}} d(z_i, \cdot) \cdot  d(z_{i+1},\cdot)] (\xi_i, \xi_{i+1}) \ .
\end{equation}
It is remains to show that the duality relation holds for the two boundary terms. i.e. that
\begin{equation}
 	[\mathcal{L}_{{\rm left}}^{{\rm BEP}}  {\mathfrak D}_T(\cdot, \xi) ](z) = [{L}_{{\rm left}}^{{\rm abs}} {\mathfrak D}_T(z, \cdot) ](\xi) \qquad \text{and} \qquad
	 [\mathcal{L}_{{\rm right}}^{{\rm BEP}}  {\mathfrak D}_T(\cdot, \xi) ](z) = [{L}_{{\rm right}}^{{\rm abs}} {\mathfrak D}_T(z, \cdot) ](\xi) \ .
\end{equation}
Being the two relations completely analogous, it is  sufficient to prove one of them, we prove it for the left boundary.
 We note that $\mathcal{L}_{{\rm left}}^{{\rm BEP}} $ acts only on site one  whereas $ [{L}_{{\rm left}}^{{\rm abs}}$ acts only on sites 0 and 1. For this reason it is sufficient to show that, for $d_\ell(k):= (T_\ell - T )^k$, 
 \begin{equation}
 	[\mathcal{L}_{{\rm left}}^{{\rm BEP}}  d_\ell(\xi_0) d(\cdot, \xi_1)](z_1) = [{L}_{{\rm left}}^{{\rm abs}} d_\ell(\cdot) d(z_1, \cdot)](\xi_0,\xi_1)
\end{equation}
 At this aim,  using the hypergeometric relation satisfied by Laguerre polynomials (see Section 9.12 in \cite{Koekoe}), we find that
\begin{align} \label{1}
z_1 \partial^{2}_{z_{1}}d(z_1, \xi_1) + \alpha  \partial_{z_{1}}d(z_1, \xi_1) = \xi_1 d(z_1, \xi_{1}-1)   \\
z_1  \partial_{z_{1}}d(z_1, \xi_1) = \xi_1 d(z_1, \xi_{1}) + \xi _1T d(z_1, \xi_{1}-1) \;. \label{2}
\end{align} 
The above identities allow us to write the action of $\mathcal{L}_{{\rm left}}^{{\rm BEP}}$ on $d(z_1, \xi_1) $ as an action on the variable $\xi_1$.
\begin{align*}
[\mathcal{L}_{{\rm left}}^{{\rm BEP}}  d_\ell(\xi_0) d(\cdot, \xi_1) ](z_1) = 
 \left(  T_\ell - T \right)^{\xi_0}  \left[  T_\ell  \xi_1 d\left(  z_1,  \xi_1 -1\right)   -  \xi_1 d\left( z_1, \xi_{1} \right)  - \xi _1 Td \left( z_1, \xi_{1}-1  \right)   \right] = \\
 \xi_1 \left[  \left(  T_\ell  - T \right)^{\xi_0 +1} d(z_1, \xi_1 -1) -  \left( T_\ell - T \right)^{\xi_0 } d(z_1, \xi_1 ) \right] =
[{L}_{{\rm left}}^{{\rm abs}} d_\ell(\cdot) d(z_1, \cdot)](\xi_0,\xi_1)
\end{align*}
that concludes the proof.
\begin{flushright}
    $\square$
\end{flushright}

\begin{remark} The so called orthogonal duality function ${\mathfrak D}_T$ is related to  the so-called generalized Laguerre polynomial via  a  normalizing factor only depending on the variable $\xi$.   More precisely, the generalized Laguerre polynomial of degree $n$, variable $x$ and parameter $ \beta$ is defined as follows
\begin{equation}
{\mathfrak L}_{\xi}^{\left( \alpha -1 \right)}(z) =\dfrac{\Gamma(\alpha + \xi) }{\Gamma(\alpha ) \xi!} \mathstrut_1 F_1 \left( {\left. \genfrac{}{}{0pt}{} {- \xi} { \alpha }  \right\vert {z}} \right)\;
\end{equation}
and then the single site duality function $d$ is related to these via the following relation
\begin{equation}
d(\zeta, k) =  (-T)^{k} \cdot \dfrac{\Gamma(\alpha ) k!} {\Gamma(\alpha +k)}\cdot {\mathfrak L}_{k}^{\left( \alpha -1 \right)}(\zeta) \ .
\end{equation}
\end{remark}

\subsubsection{Duality properties for the asymmetric process}

Once the duality relation for the symmetric process is proven we can invoke  Theorem \ref{General_result} to extend the result to the ABEP.

\begin{theorem}[Duality between open ABEP and  SIP with absorbing boundaries]
The ABEP$(\sigma,\alpha)$ with an open  boundaries, with generator ${\cal L}^{{\rm ABEP}}$ defined in  \eqref{abep}-\eqref{abep-res-dx}, is dual to the SIP$(\alpha)$ with absorbing boundaries  definded in Definition \ref{def-sip} with respect to the following duality functions:
\begin{enumerate}
\item {\bf classical duality:} 
\begin{equation} \label{dfs}
D^\sigma(x, \xi) = T_\ell ^{\xi_0} \cdot \prod_{i=1}^{N} \dfrac{\Gamma(\alpha )}{ \Gamma(\alpha  + \xi_i)} (g_{i}(x))^{\xi_i} \cdot  T_r^{\xi_{N+1}}, 
\end{equation}
\item {\bf orthogonal duality:}
\begin{equation}\label{odbepsip}
{\mathfrak D}^\sigma_T(x,\xi) = \left(  T_\ell - T \right)^{\xi_0} \cdot \prod_{i=1}^{N} (-T)^{\xi_i}\cdot \mathstrut_1 F_1 \left( {\left. \genfrac{}{}{0pt}{} {-\xi_{i}} { \alpha }  \right\vert {\frac{g_{i}(x)}{T}}} \right)\cdot
 (T_r-T)^{\xi_{N+1}}\;,
\end{equation}
for all $T >0$. Here $g $ is the map given in Definition \ref{mapg}.
\end{enumerate}
\end{theorem}
\textbf{Proof:} the result is a natural consequence of Theorem \ref{symduality} and the second item of Theorem \ref{General_result}. 
\begin{flushright}
    $\square$
\end{flushright}

\section{Applications of duality } \label{sec6}

 Due to irreducibility, the ABEP admits a unique stationary probability measure, that we will also call   steady state and we will denote it by $\mu_{ss}$. When $T_\ell = T_r=T$ this is reversible and coincides with the measure $\mu_T$ computed in Theorem \ref{rev}. When $T_\ell \neq T_r$,  reversibility is lost  and  $\mu_{ss}$   is no longer easy to compute.
We will take  advantage of the duality property proven in the previous section to compute some particular observables of $\mu_{ss}$, and more precisely, the one and two-point correlations, with respect to $\mu_{ss}$, of the observables $\{e^{-\sigma E_i(x)}, \ i\in V\}$ that are inherently related to the non-local map $g$. We will informally call these quantities $\sigma$-exponential moments or correlations. The idea is to exploit the simplicity of the dual process that  is symmetric interacting particle system. Moreover, the fact that dual particles are eventually absorbed at the boundaries, allow to compute the $\sigma$-exponential moments and correlations in terms of the absorption probabilities of the SIP particles.
\vskip.2cm
\noindent
To prove our results we use the fact that duality between two Markov generators implies duality in terms of semigroups.  This means that, if  $\{X_t\}_{t\ge0}$ and $\{Y_t\}_{t\ge0}$ are   two  Markov processes with state spaces  $\Omega$ and $\Omega^{{\rm dual}}$ respectively, whose generators are dual in the sense of Definition \eqref{gen-daulity-abcd} with respect to the duality function   $D: \Omega\times \Omega^{{\rm dual}} \to \mathbb{R}$, then for all $x\in\Omega, y\in \Omega^{{\rm dual}}$ and $t>0$,
\begin{equation}
\label{standarddualityrelation1}
\mathbb{E}_x  \big[D(X_t, y)\big]={\mathbf{E}}_{y} \big[D(x, Y_t)\big]\;
\end{equation}
where $\mathbb{E}_x $ is the expectation with respect to the
 law of the $\{X_t\}_{t\ge0}$ process started at $x$, while ${\mathbf{E}}_{y} $
denotes  expectation with respect to the law of the $\{Y_t\}_{t\ge0}$ process
initialized at $y$.

\begin{proposition}\label{p1}
Let $\mu_{ss}$ be the stationary measure of  ABEP($\sigma$, $\alpha$) with  open boundaries defined in \eqref{abep}-\eqref{abep-res-dx}, then
\begin{equation}
    \mathbb{E}_{\mu_{ss}} \left[ e^{- \sigma    E_m (x)}\right] = 1 - \sigma  \alpha \, T_\ell (N-m+1) + \dfrac{\sigma  \alpha}
    {N+1} (T_r - T_\ell) \, \dfrac{(m+N)(m-N-1)}{2}\ .
\end{equation}
\end{proposition}
\textbf{Proof:} Let $\delta_i\in \Omega^{{\rm dual}}$ the SIP($\alpha $) configuration with just one particle at site $i\in V$, then 
\begin{equation*}
    D^{\sigma}(x, \delta_i) = \dfrac{\Gamma(\alpha )}{\Gamma(\alpha +1)} \cdot  g_i(x) = \dfrac{e^{-\sigma    E_{i+1}(x)} - e^{-\sigma    E_{i}(x)} }{\sigma  \alpha} = \dfrac{e^{-\sigma    E_{i+1}(x)} }{\sigma  \alpha} \ (1-e^{-\sigma    x_i})\label{OO}
\end{equation*}
 If we initialize the dual SIP$(\alpha)$ with one particle at site $i \in V$, the dynamics can be described by a continuous time random walk $\{i(t), \ t\ge 0\}$ moving on the lattice $V \cup V^{{\rm res}}$ performing  nearest-neighbor jumps at rate $\alpha$ and absorbed at boundary sites $0$ and $N+1$. We will denote by $  \mathbb{P}_i $ the probability distribution of this process initialized at time 0 from site $i\in V$.
Then the stationary expectation of the quantity in the right hand side of \eqref{OO} linearly interpolates between $T_\ell$ and $T_r$:
\begin{align} \label{questa}
\mathbb{E}_{\mu_{ss}} \left[ e^{-\sigma    E_{i+1}(x)}  (1-e^{-\sigma    x_i}) \right] & \nn = 
      \mathbb{E}_{\mu_{ss}}  [D^{\sigma}(x, \delta_i)] = \lim_{t \to \infty} \mathbb{P}_i (i_t =0)D^{\sigma}(x, \delta_0) +   \mathbb{P}_i (i_t  =N+1)D^{\sigma}(x, \delta_{N+1})\\ & =  \sigma  \alpha \left( T_\ell + (T_r - T_\ell )\frac{i}{N+1} \right) \ .
\end{align}
    We take now the sum from $m$ to $N$ on both sides of equation \eqref{questa} to get telescopic cancellation. Since  $E_{N+1}=0$, we get
    \begin{equation}
    \mathbb{E}_{\mu_{ss}} \left[ 1-e^{- \sigma    E_m (x)}\right] = \sum_{i=m}^{N} \left( \sigma  \alpha  T_\ell + \sigma  \alpha \, (T_r - T_\ell ) \, \frac{i}{N+1}  \right)
    \end{equation}
    from which follows the result.
\begin{flushright}
    $\square$
\end{flushright}
In the next proposition we will show how to relate the above observation to gather information on the stationary  ${\sigma }$-exponential expectation of the partial energies.

\begin{remark}
Notice that the observables $\{e^{-\sigma E_i(x)}, \ i\in V\}$ are reminiscent of the microscopic Cole-Hopf transformation (known as the G{\"a}rtner transform that has been defined in \cite{Gartner} for the asymmetric exclusion process).  The Cole-Hopf transformation has been used in the literature to connect the KPZ equation for random growing interfaces and the stochastic heat equation. As remarked in \cite{Corwin}, the first hint that such transform  is available relies on the existence of a Markov duality relation.
\end{remark}

In order to compute the stationary two-point correlation of the exponential observables $\{e^{-\sigma E_i(x)}, \ i\in V\}$ we use the same strategy used in the proof of Proposition \ref{p1} to compute the $\sigma$-exponential moments. In this case, though, we initialize the dual system with two (and no longer one)  particles.

\begin{proposition}
Let $\mu_{ss}$ be the stationary measure of  ABEP($\sigma$, $\alpha$) with  open boundaries defined in \eqref{abep}-\eqref{abep-res-dx}, then
\begin{align}
    \mathbb{E}_{\mu_{ss}} \left[ e^{- \sigma     E_m (x)} e^{- \sigma    E_n (x)} \right] & = 1 - \sigma  \alpha T_\ell (2N - m- n +2) + \dfrac{\alpha  \sigma }{2(N+1)} (T_r - T_\ell ) [m^2 + n^2 -2N^2 - 2N -m-n] \nn \\ & 
    + \dfrac{(\sigma  \alpha)^2 (1-m+N)(1-n+N)}{2(N+1)(1+\alpha (N+1))}
    \left[
T_\ell^2 (N-m+2)(1+\tfrac \alpha 2  (N-n+2)) 
\right. \nn \\ & \left. + T_r^2 (N+n) (1+\tfrac \alpha 2  (N+m)) + T_\ell T_r (m(1-\alpha (n-1)) -n +\alpha (n+N(N+2)))
    \right]  \nn \\ &
   + \dfrac{(2\sigma   )^2 \alpha  (1-n+N)}{2(N+1)(1+\alpha (N+1))}    \left[
T_\ell^2 \left( \dfrac{\alpha }{3} (2n^2 + 2N^2 +2nN -n +N)  \right. \right. \nn
\\ & \left. \left.
- (n+N) [2\alpha (N+1) +1]+ 2N+1 +2\alpha (N+1)^2 \right) \right.  \nn
\\ & \left. +
T_r^2 \left( \dfrac{\alpha }{3} (2n^2 + 2N^2 +2nN -n +N) + (n+N) -1  \right) \right. \nn
\\ & \left. +
2T_\ell T_r  \left(- \dfrac{\alpha }{3} (2n^2 + 2N^2 +2nN -n +N) + (n+N)(\alpha (N+1) -1) +1  \right)
\right] \nn
\end{align}
where $m \leq n$.
\end{proposition}
\textbf{Proof:}
Let $\xi= \delta_i + \delta_j\in \Omega^{{\rm dual}}$ be the dual configuration with two particles at sites $i,j\in V$, $i \neq j$. The duality function evaluated in $\xi$ is then given by 
\begin{equation}
    D^{\sigma} (x, \delta_i + \delta_j) = \dfrac{e^{-\sigma    E_{i+1}(x)} - e^{-\sigma    E_{i}(x)} }{\sigma  \alpha} \cdot \dfrac{e^{-\sigma    E_{j+1}(x)} - e^{-\sigma    E_{j}(x)} } {\sigma  \alpha} \ .
\end{equation}
Considering the expectation with respect to the stationary measure:
\begin{align}  \label{2punti}
& \mathbb{E}_{\mu_{ss}} \left[  \left( e^{-\sigma    E_{i+1}(x)} - e^{-\sigma    E_{i}(x)} \right) \left(e^{-\sigma    E_{j+1}(x)} - e^{-\sigma    E_{j}(x) } \right) \right] = 
      ( \sigma  \alpha)^2 \cdot   \mathbb{E}_{\mu_{ss}}[D^{\sigma}(x, \delta_i + \delta_j) ]\\ & \nn 
      = ( \sigma  \alpha)^2 \lim_{t \to \infty} \bigg\{\mathbb{P}_{i,j} (i_t =0, j_t = 0 )D^{\sigma}(x, 2\delta_0) +   \mathbb{P}_{i,j} (i_t  =N+1, j_{t} = N+1)D^{\sigma}(x, 2\delta_{N+1}) +   \\ \nn &  D^{\sigma}(x, \delta_0 + \delta_{N+1}) \left( \mathbb{P}_{i,j} (i_t  =0, j_{t} = N+1) + \mathbb{P}_{i,j} (i_t  =N+1, j_{t} = 0)
      \right)\bigg\}
      \\ & = ( \sigma  \alpha)^2 \left\{ T_\ell^2 \dfrac{[1+\alpha (N+1-i)](N+1-j)}{(N+1)(1+\alpha (N+1))} + 
      T_r^2 \dfrac{i(1+\alpha j)}{(N+1)(1+\alpha (N+1))}\right. \\&+\left.
      T_\ell T_r  \dfrac{[\alpha (N+1) -1]i +[1+\alpha (N+1)]j - 2\alpha ij}{(N+1)(1+\alpha (N+1))}
      \right\}
\end{align}
where $\mathbb{P}_{i,j} $ is the probability distribution associated to two dual SIP$(\alpha)$ particles $\{(i(t),j(t)), \ t\ge 0\}$.
On the other hand, if $i=j$ we have:
\begin{equation}
    D^{\sigma} (x, 2\delta_i ) = \dfrac{ \left(e^{-\sigma    E_{i+1}(x)} - e^{-\sigma    E_{i}(x)}  \right)^2}{\alpha  (\alpha +1) \sigma ^2} 
\end{equation}
and considering the expectation with respect to the stationary measure:
\begin{align}  \label{2punti}
& \mathbb{E}_{\mu_{ss}} \left[  \left( e^{-\sigma    E_{i+1}(x)} - e^{-\sigma    E_{i}(x)} \right)^2 \right] = 
      \mathbb{E}_{\mu_{ss}}  \alpha  (\alpha +1) \sigma ^2 D^{\sigma}(x, 2\delta_i ) \\ & \nn 
      = \alpha  (\alpha +1) \sigma ^2 \lim_{t \to \infty} \bigg\{\mathbb{P}_{i,i} (i_t =0, i_t = 0 )D^{\sigma}(x, 2\delta_0) +   \mathbb{P}_{i,i} (i_t  =N+1, i_{t} = N+1)D^{\sigma}(x, 2\delta_{N+1}) +   \\ \nn &  D^{\sigma}(x, \delta_0 + \delta_{N+1}) \left( \mathbb{P}_{i,i} (i_t  =0, i_{t} = N+1) + \mathbb{P}_{i,i} (i_t  =N+1, i_{t} = 0)
      \right)\bigg\}
      \\ & =  \alpha  (\alpha +1) \sigma ^2 \left\{ T_\ell^2 \frac{2(N+1-i) (\alpha  (N+1-i)+1) -1}{2(N+1) (\alpha  (N+1)+1)}  + \right. \\& \left.
      T_r^2 \dfrac{2i(1+\alpha i) -1}{2(N+1)(\alpha (N+1)+1)}  
      + T_\ell T_r  \dfrac{(\alpha (N+1)-1)i + (\alpha (N+1)-1)i -2\alpha i^2 +1}{(N+1)(\alpha (N+1)+1)} 
      \right\}
      \nonumber 
          \\&=
     \alpha  (\alpha +1) \sigma ^2 \left\{ T_\ell^2 \frac{ 2\alpha  i^{2}+(-4\alpha  N-4\alpha +2)i+(2\alpha  N^2+4\alpha  N+2\alpha -2 N-3)}{2(N+1) (\alpha  (N+1)+1)}  + \right. \\& \left.
      T_r^2 \dfrac{2\alpha  i^2+2 i-1}{2(N+1)(\alpha (N+1)+1)}  
      + T_\ell T_r  \dfrac{-2\alpha  i^2 + 2i(\alpha  N+\alpha -1) +  1}{(N+1)(\alpha (N+1)+1)}
      \right\}
      \nonumber
\end{align}
This allows us to gather informations on the two-point $\sigma$-exponential stationary correlations. To achieve this 
we take a double sum in equation \eqref{2punti},  one from $m$ to $N$ and one from $n$ to $N$. By telescopic arguments one then gets 
\begin{align}\label{ultima}
&  \mathbb{E}_{\mu_{ss}} \left[e^{- \sigma     E_m (x)} e^{- \sigma    E_n (x)} \right] = \mathbb{E}_{\mu_{ss}} \left[e^{- \sigma     E_m (x)}\right] + \mathbb{E}_{\mu_{ss}} \left[e^{- \sigma     E_n (x)}\right]  - 1 
 +\nn \\ & \nn (\sigma  \alpha)^2 \sum_{i=m}^{N} \sum_{j=n}^{N} \left\{ T_\ell^2 \mathbb{P}_{i,j} (i_t  =0, j_{t} = 0) + T_r^2 \mathbb{P}_{i,j} (i_t  =N+1, j_{t} = N+1) \right.\\&+\left. T_\ell T_r \left[ \mathbb{P}_{i,j} (i_t  =0, j_{t} = N+1)+ \mathbb{P}_{i,j} (i_t  =N+1, j_{t} = 0) \right] \right\} \nn \\ &
 + (2\sigma   )^2 \alpha  \sum_{i=n}^{N} \left\{ T_\ell^2 \mathbb{P}_{i,i} (i_t  =0, i_{t} = 0) + T_r^2 \mathbb{P}_{i,i} (i_t  =N+1, i_{t} = N+1) \right. \nn\\&+\left. T_\ell T_r \left[ \mathbb{P}_{i,i} (i_t  =0, i_{t} = N+1)+ \mathbb{P}_{i,i} (i_t  =N+1, i_{t} = 0) \right] \right\} 
 \end{align}
 where the first two terms on the right hand side have been computed in the previous theorem.
 To conclude the proof it remains to  are plug in the expression above  the absorption probabilities of two dual SIP$\alpha$ particles absorbed at the boundaries 0 and $N+1$. 
 These are harmonic function of the two dimensional Laplacian. They solve a systems of discrete equations with appropriate boundary conditions. We show how to get $p_{i,j}:=  \mathbb{P}_{i,j} (i_t  =0, j_{t} = 0)$  for $i,j \in V $ as the others can be found similarly.
 \begin{equation}
     \begin{cases}
         4p_{i,j} = p_{i-1,j} + p_{i+1,j} + p_{i,j-1} + p_{i,j+1} \\
         2 p_{i,i} = p_{i-1,i} + p_{i,i+1} 
     \end{cases}
 \end{equation}
 for the first two equations we get that
 \begin{equation*}
    p_{i,j}= Ai + Bj + Cij + D \quad \text{for} \quad i \neq j 
 \end{equation*}
 and 
  \begin{equation*}
    p_{i,i}= (A+B)i + Ci^2 + D + \frac{B-A}{2} \ .
 \end{equation*}
Three of the unknown can be found using the boundary conditions:
  \begin{equation}
     \begin{cases}
        p_{0,0} = D = 1 \\
        p_{0,j} = Bj + D = 1- \frac{j}{N+1} \\
          p_{N+1,N+1} = A(N+1) + B(N+1) + C(N+1)^2 + D = 0
     \end{cases}
 \end{equation}
 while the last one can be found conditioning on the first jump, i.e. 
 \begin{equation*}
     (4\alpha  +2)   p_{i,i+1} = \alpha    p_{i-1,i+1} +   \alpha  p_{i,i+2} + (\alpha +1)   p_{i,i} +   (\alpha  + 1 )p_{i+1,i+1} \;.
 \end{equation*}
This leads to the following solutions for the four unknown  
\begin{equation}
     \begin{cases}
     A=-\frac{\alpha }{1+\alpha (N+1)}\\
      B = - \frac{1}{N+1} \\
       C=\frac{\alpha }{(1+N)(1+\alpha (N+1)} \\
        D = 1   \end{cases}
\end{equation}
Finally we obtain 
\begin{equation}
     \mathbb{P}_{i,j} (i_t  =0, j_{t} = 0) = p_{i,j}=\frac{(N+1-j) (\alpha  (-i+N+1)+1)}{(N+1) (\alpha  (N+1)+1)} - \dfrac{1}{2(N+1)(\alpha (N+1)+1)} \mathbbm{1}_{\{ i=j \} }
\end{equation}
for the absorption probabilities of both particles to the left. Similarly one can get the absorption probabilities of both particles to the right: 
\begin{equation}
    \mathbb{P}_{i,j} (i_t  = N+1, j_{t} = N+1) = p_{i,j}=\dfrac{i(1+\alpha j)}{(N+1)(\alpha (N+1)+1)}  - \dfrac{1}{2(N+1)(\alpha (N+1)+1)} \mathbbm{1}_{\{ i=j \} }
\end{equation}
and the absorption probability of one particle to the left and one to the right
\begin{equation}
\begin{split}
    & \mathbb{P}_{i,j} (i_t  =0, j_{t} = N+1) +  \mathbb{P}_{i,j} (i_t  =N+1, j_{t} = 0)\\=& p_{i,j}=\dfrac{(\alpha (N+1)-1)i + (\alpha (N+1)-1)j -2\alpha ij}{(N+1)(\alpha (N+1)+1)}  + \dfrac{1}{(N+1)(\alpha (N+1)+1)} \mathbbm{1}_{\{ i=j \} }\ . 
     \end{split}
\end{equation}
Substituting these expressions in \eqref{ultima} we obtain the result. 
\begin{flushright}
    $\square$
\end{flushright}

\subsection*{Acknowledgments}
The authors benefited from inspiring conversation with Leonardo de Carlo, Cristian Giardin\`{a}, Claudio Giberti, Seth Lloyd and Frank Redig.
C.F. acknowledges hospitality and support from the Galileo Galilei Institute (GGI) during the scientific program `Randomness, Integrability and Universality' held in Spring 2022.

\subsection*{Data Avaibility Statement}
No new data were created or analysed in this study.

\end{document}